\documentclass[12pt]{article}
\usepackage{amsfonts,amssymb,amsmath,graphicx,float,mathtools,amsthm}

  \newcommand{\lcr}{\raisebox{-5pt}{\mbox{}\hspace{1pt}
                 \includegraphics{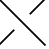}\hspace{1pt}\mbox{}}}
\newcommand{\ift}{\raisebox{-5pt}{\mbox{}\hspace{1pt}
                 \includegraphics{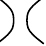}\hspace{1pt}\mbox{}}}
\newcommand{\zer}{\raisebox{-5pt}{\mbox{}\hspace{1pt}
                 \includegraphics{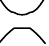}\hspace{1pt}\mbox{}}}

    \newtheorem{theorem}{Theorem}[section]
    \newtheorem{prop}[theorem]{Proposition}
    \newtheorem{lemma}[theorem]{Lemma}
    
    \newtheorem{rem}[theorem]{Remark}
    \newtheorem{definition}[theorem]{Definition}
    
\newtheorem{example}[theorem]{Example}

\title{Locality and the Uniqueness of Quantum Invariants}
\author{Charles Frohman \\ Department of Mathematics, The University of Iowa}

\begin{document}

\maketitle

\begin{abstract}We introduce the notion of a "state function" for framed tangles in a disk. After choosing a finite set of states for each marked disk, a state function is a projection from the vector space spanned by all tangles to the vector space spanned by the states, that is local, and topologically invariant. Given the states for the Kauffman bracket, and the quantum $SU(3)$-invariant we classify all state functions, and then compare our results to the literature.\end{abstract}

\section{Introduction} 

In linear algebra a function from the Cartesian product of vector spaces to the base field is called {\bf multilinear} if fixing all but one variable gives a linear functional. Analogously, given links that agree outside a ball, a linear relation between their invariants is called a {\bf skein relation}.  In the mid-eighties Lickorish \cite{L}  used the idea of the {\bf skein} of a planar region to codify multilinearity in the theory of knot polynomials. Kauffman understood that skein relations were associated to  partition functions from quantum field theory, and used them both for the diagrammatic study of knots and links, and to make topological arguments about quantum field theory \cite{K1}.

In quantum field theory a partition function is {\bf local} if you can compute its global value by combining values computed in small pieces of space-time.
   If a partition function is topologically invariant and local, its values satisfy skein relations.  Witten  used locality to realize the  Jones polynomial as a partition function in a topological quantum field theory \cite{W}.  The heart of his argument (page 378) is  that up to a scalar multiple, there is a unique one parameter family  of  invariants of framed links based on the data from the representation theory of the Lie algebra $sl_2$. The invariants only exist for certain values of the parameter, however they  satisfy the skein relation for the Jones polynomial. We present an elementary descendant of that argument here, avoiding physical arguments and representation theory. The uniqueness theorems in this paper are similar to the axiomatic characterization of the determinant \cite{HK}.

At the end of this paper  we use the uniqueness theorem to make sense of the plethora of definitions of quantum $SU(3)$-invariants. The definition of a quantum invariant can have different versions depending on how vertices and crossings in diagrams are tracked.  For instance if the link is oriented it is common to multiply  the invariant by some function of the total writhe of the link in order to get an invariant of oriented links as opposed to framed links. A similar flexibility in definition is seen when diagrams are used to represent invariant tensors. If a tensor is invariant then any scalar multiple of it is invariant, and this allows scaling by a function that is a product of contributions from the vertices. In an operation called {\bf twisting}, equivalent formulations result although the contribution of a vertex can depend on how it has been decorated.  {\em In the two examples worked in this paper, the invariants of tangles, up to twisting, depend on a single variable.}  While there are  more robust frameworks for expressing locality,  working in such a framework would obscure the central idea of this paper, which is that invariants can be obtained by solving relatively simple systems of equations arising from just locality and the Reidemeister moves. The state functions are defined by giving a collection of reduction rules that are confluent and terminal \cite{SW}.

The representation theoretic approach to the construction of quantum invariants of framed links involves a {\bf quantum group} \cite{RT} , which is a deformation of the universal enveloping algebra of a Lie algebra. Over the ring of formal power series the deformation is trivial, so the category of representations of the quantum group is isomorphic to the category of representations of the Lie algebra. As a result the information about quantum invariants should be carried by the original Lie algebra.

The connection between diagrams in a disk and representation theory was codified by Kuperberg \cite{Ku}, with the introduction of the concept of a {\bf spider}. Kuperberg gave presentations of the spiders associated to rank $2$ Lie algebras.  A presentation of the spider for $A_n$ is given in \cite{CKM}. Modeling the invariant tensors of the Lie algebra via a spider, reduces the study of quantum invariants to the study of polynomial equations derived from Reidemeister moves on diagrams.  
The uniqueness of a local invariant given the  data from the representation theory of a Lie algebra was used  implicitly in \cite{Ku,MOY} and \cite{OY} to determine link polynomials associated with $SU(n)$, $A_2$ and $G_2$. The idea is further explored in \cite{MPS},  along with an exhaustive analysis of categories of diagrams with a single trivalent vertex.

This paper is an elementary introduction to a program for  constructing  skein algebras associated to arbitrary Lie algebras and super Lie algebras from a diagrammatic viewpoint.  In a less elementary follow up, we will abstract the notion of state function here, and apply it to the construction of skein algebras, and the explication of their structure.
 
 The second section begins by reviewing the diagrammatic presentation of framed links, and the Kauffman bracket.  It continues by giving the definition of a {\em state function} which is a combinatorial analog of a partition function from quantum field theory.  As a warm up, we cast Kauffman's uniqueness argument in the light of state functions, proving it is the unique separating state function.

 In the last section we prove a  uniqueness theorem for quantum $SU(3)$-invariants of framed knotted trivalent graphs, and use it to compare the different versions that appear in the literature \cite{Ku,OY,Kh,Si}. 
 
 I would like to thank Fred Goodman, Joanna Kania-Bartoszynska, Lou Kauffman and Shawn Nevalainen for helpful input.

\section{The Kauffman bracket}

\subsection{Diagrams of Framed Links}

A link in $\mathbb{R}^3$ is encoded by a four valent planar graph, with data at the vertices indicating over-crossings.  The link can be reconstructed up to isotopy from such a diagram. 
Two diagrams represent the same link if you can pass from one to the other by isotopies that preserve the combinatorial data ( sometimes called R0 moves) and three different moves that change the combinatorics of the diagram \cite{CF}. 
These are the Reidemeister moves.  The first Reidemeister move, or RI is
\begin{equation} \includegraphics{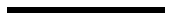} \longleftrightarrow \raisebox{-15pt}{ \includegraphics{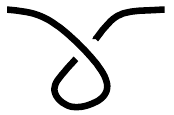}}. \end{equation}
The second Reidemeister move or RII is,
\begin{equation}  \raisebox{-15pt}{\includegraphics{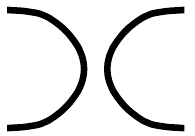}} \longleftrightarrow \raisebox{-15pt}{ \includegraphics{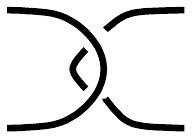}}. \end{equation}
Finally, the third Reidemeister move or RIII is 
\begin{equation}  \raisebox{-15pt}{\includegraphics{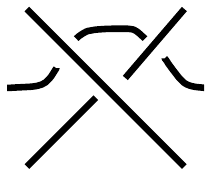}} \longleftrightarrow \raisebox{-15pt}{ \includegraphics{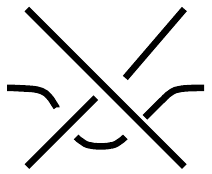}}. \end{equation}
 A {\bf framed link} in $\mathbb{R}^3$ is an
embedding of a disjoint union of annuli into $\mathbb{R}^3$.    Diagrammatically we  depict framed links by showing the core of the annuli. You should imagine the annuli  lying parallel to the plane of the paper, this is sometimes called {\bf the blackboard framing}. Two framed links in $\mathbb{R}^3$ are equivalent if they are isotopic.

Every framed link can be isotoped so that it has a projection into the plane that is blackboard framed.   Two diagrams represent the same framed link with the blackboard framing if you can pass from one to the other by isotopies that don't change the combinatorics of intersection, the second Reidemeister mover RII, and the third Reidemeister move RIII \cite{K2}.

\subsection{State sum description of the Kauffman bracket}

There are two ways of smoothing a crossing of a diagram.

\begin{figure}[H]\begin{center} \includegraphics{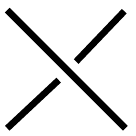} \hspace{.5in}  \includegraphics{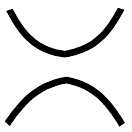} \hspace{.5in} \includegraphics{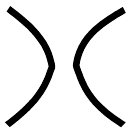} \caption{Smoothing a Crossing}\end{center} \label{smooth}\end{figure}

The middle diagram in Figure $1$ shows the positive smoothing of the crossing, the right diagram shows the negative smoothing. Here is how to tell then apart. Put your hand on the over-crossing arc in the diagram, and sweep counterclockwise. The pair of opposite corners that you sweep out, are the ones that are joined in the positive smoothing.

If the diagram has $n$ crossings, then there are $2^n$ ways of smoothing the diagram.  These are the {\bf states} of the diagram.  Given a state $s$, there are three obvious numbers.  The first, $p(s)$, is the number of positive smoothings that were made to arrive at the state. The second, $n(s)$, is the number of negative smoothings. Notice that $p(s)+n(s)$ is the number of crossings of the diagram.  The third number, $c(s)$, is the number of connected components of the state.

Suppose that $D$ is a diagram and $S$ is the set of states of the diagram. The {\bf  Kauffman bracket} of the diagram $D$ is given by 
\begin{equation} <D>=\sum_{s\in S} A^{p(s)-n(s)}(-A^2-A^{-2})^{c(s)},\end{equation}
where the sum takes place  in the ring of Laurent polynomials with integer coefficients in the variable $A$, $\mathbb{Z}[A,A^{-1}]$.

The Kauffman bracket is an invariant of the framed link that $D$ represents \cite{K}.  Kauffman proved the invariance by showing that the value is unchanged by the Reidemeister moves for framed links.  The first step was  to derive {\bf skein relations} that the Kauffman bracket satisfies.  A skein relation is a linear relation between diagrams that agree outside of a disk, with coefficients in some ring. In this case the ring is $\mathbb{Z}[A,A^{-1}]$.  We depict the relation by only showing their intersection in the disk where the diagrams are different. 

\begin{theorem} \cite{K} \label{skeink} The Kauffman bracket satisfies the skein relations
  \[\left<\lcr\right>-A\left<\zer\right>-A^{-1}\left<\ift\right>=0 \]
  and \begin{equation*}\left<\bigcirc \cup L\right>+(A^2+A^{-2})\left<L\right>=0.\end{equation*} \end{theorem} \qed

Actually Kauffman defined the bracket of a diagram $D$ with states $S$ to be
\begin{equation} <D>=\sum_{s\in S} A^{p(s)}B^{n(s)}d^{c(s)-1},\end{equation}
and then showed that it could only be invariant under the second Reidemeister move if $B=A^{-1}$ and $d=-A^2-A^{-2}$. He chose the exponent $c(s)-1$ for $d$, so that his final answer to match the normalization of the Jones polynomial. 

\subsection{State Functions and Uniqueness}

In this section we work over the complex numbers. Hence
\begin{equation} A \in \mathbb{C}-\{0\}.\end{equation}

Let $(D^2,\mathbf{v})$ be an oriented disk $D^2$ , decorated with points $\mathbf{v}=\{v_1,v_2,\ldots,v_n\}\subset \partial D^2$.  We begin by considering collections of disjoint, proper embeddings of arcs, that are transverse to $\partial D$, and whose boundary is the set ${\bf{v}}$.  A {\bf crossingless matching} is an isotopy class, relative to $\partial D$ of a choice of such arcs. A crossingless matching is determined by which pairs of points in $\mathbf{v}$ make up the boundaries of  a set of arcs that represent the crossingless matching.  Since each arc has two endpoints, unless the cardinality of ${\bf{v}}$ is even, the set of crossingless matchings is empty. The number of crossingless matchings with $2n$ boundary points is equal to the $n$th {\bf Catalan number},
\begin{equation} \frac{1}{n+1}\binom{2n}{n}.\end{equation}

The state space  $\mathcal{S}(D,\mathbf{v})$ is the complex vector space  having basis the crossingless matchings.  The state space of an oriented disk with no points on the boundary  is $\mathbb{C}\emptyset$, where $\emptyset$ is the empty diagram.

Next we describe  the set of {\bf tangle  diagrams} in the disk with boundary points ${\mathbf{v}}$. These are graphs with four-valent and mono-valent vertices, embedded in $D$ where the mono-valent vertices are exactly the set $\mathbf{v}$, and the four-valent vertices have over-crossing data. We assume the arc and circle components of a tangle diagram are smooth and the arc components are transverse to $\partial D$.  Two tangle diagrams are equivalent if you can pass from one to the other by isotopies relative to $\partial D$ that preserve the crossing data, and the moves RII, and RIII. We call the equivalence classes {\bf framed homogeneous tangles}. Let    $\mathcal{F}(D,{\bf v})$ denote the complex vector space with basis
 the framed homogeneous tangles with boundary ${\bf{v}}$. Elements of  $\mathcal{F}(D,{\bf v})$ are treated as linear combinations of diagrams, up to isotopy and the second and third Reidemeister moves.  Framed homogeneous tangles correspond to isotopy classes of  embeddings of annuli and strips into $D\times [0,1]$, so that the same side of each strip is up at its two ends.
 
 Suppose that $F:D\rightarrow D'$ is a smooth orientation preserving diffeomorphism of disks.  Notice that $F$ acts on crossingless matchings, and extending linearly it gives rise to 
 \begin{equation} F:S(D,\mathbf{v})\rightarrow S(D',F(\mathbf{v})).\end{equation}

 The diffeomorphism $F$ can be can be extended to an orientation preserving diffeomorphism $\overline{F}:D\times [0,1]\rightarrow D'\times [0,1]$ by letting $\overline{F}(x,t)=(F(x),t)$. 
 Similarly $\overline{F}$ extends linearly to yield
  \begin{equation} \overline{F}:\mathcal{F}(D,\mathbf{v})\rightarrow \mathcal{F}(D',F(\mathbf{v})).\end{equation}  The action of $\overline{F}$ on tangles is determined by the action of $F$ on tangle diagrams.

Suppose that  for each oriented disk $D$ and $\mathbf{v}\subset \partial D$ we have chosen a linear map
\begin{equation} Z_D:\mathcal{F}(D,\mathbf{v})\rightarrow \mathcal{S} (D,{\mathbf{v}}).\end{equation} 

We say the set of choices $\{Z_D\}$ is {\bf topologically invariant} if for all orientation preserving diffeomorphisms $F:D\rightarrow D'$, and all framed homogeneous tangles $T$
\begin{equation} F\circ Z_{D}(T)=Z_{D'}(\overline{F}(T)).\end{equation}

Let   $D'\subset D$  be a disk given the orientation from $D$ whose boundary is in {\bf general position} with respect to the tangle diagram $T$. This means that $\partial D'$ is disjoint from the vertices of $T$ and intersects the edges of $T$ transversely.  We can apply $Z_{D'}$ to that part of $T$ that lies in $D'$ and glue the result in a smooth way into the part of $T$ lying outside of $D'$. We abuse notation by denoting the result $Z_{D'}(T)$.

\begin{example}Assume
\begin{equation} Z_{D'}(\raisebox{-8pt}{\scalebox{.5}{\includegraphics{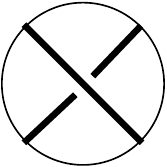}}})=A\raisebox{-8pt}{\scalebox{.5}{\includegraphics{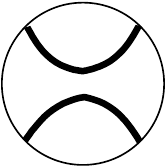}}} +A^{-1}\raisebox{-8pt}{\scalebox{.5}{\includegraphics{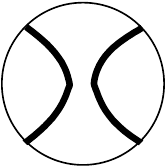}}} \end{equation}

In Equation \ref{locpic} the disk $D'$ is dotted. On the right hand side we show the result of gluing $Z_{D'}(T\cap D')$ into $T \cap(D-D')$.

\begin{equation}\label{locpic}  Z_{D'}(\raisebox{-15pt}{\scalebox{.5}{\includegraphics{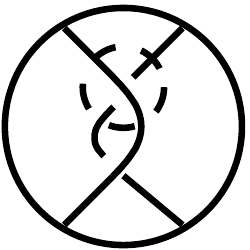}}})=A\raisebox{-15pt}{\scalebox{.5}{\includegraphics{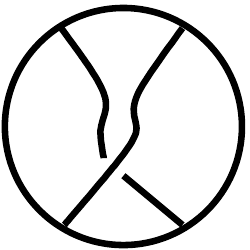}}} +A^{-1}\raisebox{-15pt}{\scalebox{.5}{\includegraphics{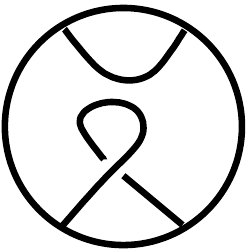}}} \end{equation} \end{example}

We say that the set of choices of $Z_D$ is {\bf local} if for all $D'\subset D$ and $T$ with $T$ transverse to $\partial D'$, \begin{equation}Z_D(T)=Z_D(Z_{D'}(T)).\end{equation}
\begin{definition} A choice of linear maps  for each oriented marked disk \begin{equation} Z_D:\mathcal{F}(D,\mathbf{v})\rightarrow \mathcal{S} (D,{\mathbf{v}}).\end{equation}  for  is a {\bf state function} if it is topologically invariant and  local. \end{definition}

The formula for the Kauffman bracket leads to  a state function. Choose $A\in \mathbb{C}-\{0\}$.  Let $T$ be a tangle diagram with boundary ${\bf{v}}\in \partial D^2$ and crossings $S$. A state $s$ is a choice of smoothings for all the crossings of $T$. Let $p(s)$ denote the number of positive smoothings, $n(s)$ denote the number of negative smoothings.  The resulting diagram  consists of some circles, and a crossingless matching $m$. Let $c(s)$ denote the number of circles. 
\begin{equation} <D>_A=\sum_{s\in S} A^{p(s)-n(s)}(-A^2-A^{-2})^{c(s)}m\in  \mathcal{S} (D,{\mathbf{v}}).\end{equation}
The linear extension of this formula yields,
\begin{equation} <\ >_A:  \mathcal{F}(D,{\mathbf{v}})\rightarrow \mathcal{S} (D,{\mathbf{v}})\end{equation} the {\bf Kauffman bracket state function}.
Notice that for a crossingless matching $m$,
$<m>_{A}=m$. The proof that it is well defined is exactly the proof that the Kauffman bracket is an invariant of framed links.

From the diagrammatic viewpoint it is clear that
\begin{equation} \mathcal{S}(D,{\bf{v}})\leq \mathcal{F}(D,{\bf{v}}). \end{equation}
crossingless matchings are just equivalence classes of tangle diagrams without crossings. Distinct crossingless matchings are not related by Reidemeister moves  because they are determined by the endpoints of the arcs making up the crossingless matching.

\begin{lemma}  \label{project} Every state function  is a  projection. \end{lemma}

\proof  Let $(D,{\bf{v}})$ be the unit  disk in $\mathbb{R}^2$ with the induced orientation and  an even number of points ${\bf{v}}$  chosen in its boundary.  Let $\mathcal{B}=\{m_1,\ldots,m_n\}$ be the ordered set of crossingless matchings with boundary ${\bf{v}}$.  They form a basis for $\mathcal{S}(D,{\bf{v}})$.  Since $\mathcal{S}(D,{\bf{v}})\leq \mathcal{F}(D,{\bf{v}})$ the restriction of $Z_D$ to  $\mathcal{S}(D,{\bf{v}})$ can be represented as a matrix with respect to $\mathcal{B}$,
\begin{equation} [Z_D]_{\mathcal{B}}.\end{equation} 

Choose the $m_i$ so that there is a disk $D_r$ of  radius $r<1$ centered at the origin  so that $\partial D_r$ is transverse to all the $m_i$ and an orientation preserving  diffeomorphism $F:D_r\rightarrow D$ so that $F(m_i\cap D_r)=m_i$. Let ${\bf{v}}_r$ be the set of points $m_1\cap \partial{D}_r$. The invariance of $Z$ under diffeomorphism implies that the matrix representing $Z_{D_r}$ restricted to $\mathcal{S}(D_r, {\bf{v}}_r)$ with respect to $\mathcal{B}_r=\{m_i\cap D_r\}$ is the same as $ [Z_D]_{\mathcal{B}}$. 
We can evaluate each $m_i$ by first evaluating inside $D_r$ and then evaluating again.  This means that
\begin{equation} [Z_D]_{\mathcal{B}}^2=[Z_D]_{\mathcal{B}}.\end{equation}  By the definition of projection, $Z_D$ is a projection onto its image. \qed
 
A state function is {\bf separating} if its restriction to $\mathcal{S}(D,\mathbf{v})$ always has full rank.

\begin{theorem}  Any separating state function of framed homogeneous tangles in a disk is equal to  the Kauffman bracket state function for some choice of $A\in\mathbb{C}-\{0\}$. \end{theorem}

\proof  Let $D$ be a disk with four points ${\bf{v}}$ in its boundary.  The state space of $(D,{\bf{v}})$ is a two dimensional vector space with basis the two crossingless matchings with four endpoints. Let $Z$ denote the state function.

Since the crossingless matchings form a basis for  $(D,{\bf{v}})$ the diagram with two components and a single crossing evaluates as
\begin{equation}\label{assume} Z\left(\raisebox{-20pt}{\includegraphics{crosst.pdf}}\right)=A \raisebox{-20pt}{ \includegraphics{zerot.pdf}} + B\ \raisebox{-20pt}{\includegraphics{infinityt.pdf}}, \end{equation}
for some $A,B\in\mathbb{C}$.
A trivial component in a disk must evaluate as a multiple of the empty diagram,
\begin{equation} \label{assume2} Z\left(\raisebox{-20pt}{\includegraphics{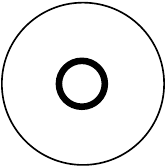}}\right)=d\emptyset.\end{equation}
Here we are denoting the empty crossingless matching by $\emptyset$.  Locality requires that $Z(\emptyset)=\emptyset$, otherwise you could evaluate inside a trivial disk and change the value of a tangle.

To be an invariant of framed homogeneous tangles, the partition function must send two tangles that differ by the second Reidemeister move to the same linear combination of  crossingless matchings. Locality means that this can be analyzed by isolating the second Reidemeister move in a disk, and noting
\begin{equation} Z\left(\raisebox{-20pt}{\includegraphics{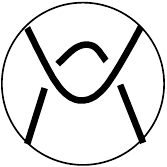}}\right)= Z\left(\raisebox{-20pt}{\includegraphics{zerot.pdf}}\right)\end{equation}
In order to expand the left hand side we isolate the two crossings in small disks and expand using Equation \ref{assume}, and use Equation \ref{assume2} to eliminate the trivial component.  Collecting terms we get 
\begin{equation} ABZ\left(\raisebox{-20pt}{\includegraphics{zerot.pdf}}\right)+(A^2+B^2+d)Z\left(\raisebox{-20pt}{\includegraphics{infinityt.pdf}}\right)= Z\left(\raisebox{-20pt}{\includegraphics{zerot.pdf}}\right).\end{equation}
Linear independence forces $AB=1$, so $B=A^{-1}$ and $d=-A^2-A^{-2}$.  An argument by induction on the number of crossings shows that  for every $(D,{\bf{v}})$, $Z_D$  the Kauffman bracket state function where the variable is $A$.   \qed

We have not addressed existence.   The skein relations in Theorem \ref{skeink} give rise to confluent, terminal reduction rules,
\begin{equation} S_1:  \left<\lcr\right>\rightarrow A\left<\zer\right>+A^{-1}\left<\ift\right>\end{equation}
and 
\begin{equation} S_2: \bigcirc\rightarrow (-A^2-A^{-2})\emptyset .\end{equation}
To compute the partition function for any disk $D$, start by placing small balls around all the crossings, and applying the first reduction rule. On each of the resulting diagrams but disks around each trivial component and apply the second reduction rule.

The Kauffman bracket {\em state function} for an oriented disk with no boundary points takes on values in $\mathbb{C}\emptyset$.  The Kauffman bracket takes on values in $\mathbb{C}$. You can normalize the Kauffman bracket by choosing a linear functional on $\mathbb{C}\empty$. This is equivalent to choosing a value of the empty link.    For the volume conjecture it makes sense to  send $\emptyset$ to $\frac{-1}{A^2+A^{-2}}$,  for the correspondence with representation theory $1$ is a better value.

\begin{rem}A version of the uniqueness theorem holds over $\mathbb{Z}[A,A^{-1}]$. However, the conclusion is that any  separating  state function is the result of evaluating the Kauffman bracket at $(-1)^kA^l$ for some $k,l\in \mathbb{Z}$.\end{rem}

\section{Quantum $SU(3)$-invariants}

\subsection{Webs and Tangles}
A {\bf  web} is an oriented, graph $\Gamma$ embedded in a disk $D$  having only trivalent and mono-valent vertices. The edges have been oriented so that each trivalent vertex is a source or sink. The mono-valent vertices coincide with $\Gamma\cap \partial D$, and the edges of $\Gamma$ intersect $\partial D$ transversely. The web decomposes the disk into vertices, edges and faces.  A face is {\bf interior} if none of its sides are contained in $\partial D$. The fact that the vertices are sources or sinks implies that the orientations of the sides of an interior face alternate as you traverse its boundary. This means that interior faces have an even number of sides. A {\bf web} is non-elliptic \cite{Ku} if its interior faces have at least six sides.

We now consider oriented disks $D$ along with a choice of decorated points ${\bf{v}}\subset
\partial D$. The decoration is a choice of sign $\pm$ for each point $v\in {\bf{v}}$.  A web $\Gamma$ in $(D,{\bf{v}})$ is a non-elliptic web in $D$ whose intersection with $\partial D$ is ${\bf{v}}$, so that if the sign of $v$ is positive the edge of $\Gamma$ ending at $v$ points into the disk, and if the sign of $v$ is negative then the edge of $\Gamma$ ending there points out of $D$.  Let $\mathcal{S}(D,{\bf{v}})$ be the complex vector space with basis the set of isotopy classes relative to $\partial D$ of non-elliptic webs in $(D,{\bf{v}})$.  We call $\mathcal{S}(D,{\bf{v}})$ the space of {\bf states} of $(D,{\bf{v}})$.

In order to capture knotting, we consider {\bf tangle diagrams} that are  oriented graphs $\Gamma$ embedded in $(D,{\bf{v}})$ that have, in addition to mono-valent and trivalent vertices,  four valent vertices that carry crossing data. The orientations of the edges at the crossings match up, so that when you push the graph into $\mathcal{D}\times [0,1]$ according to the crossing data, the edges are oriented coherently. Just as before, the vertices are expected to be sources and sinks, and the signs on $v\in {\bf{v}}$ correspond to whether the edges point in or out at $v$.

Two diagrams are  equivalent it they differ by  isotopy of the diagrams relative to $\partial D$  that don't  change the decorations or crossing data, the second and third Reidemeister moves, and the moves obtained from the moves below with orientations added, and perhaps the crossings all simultaneously changed.
\begin{equation} \label{newr} \raisebox{-15pt}{\includegraphics{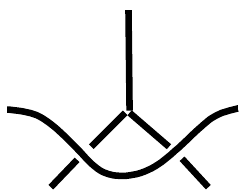}} \longleftrightarrow \raisebox{-15pt}{ \includegraphics{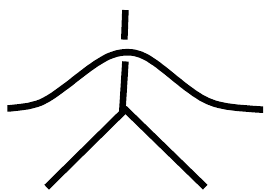}}. \end{equation}
A {\bf framed homogeneous tangle} is  an equivalence class of diagrams with respect to this notion of equivalence.

Let $\mathcal{F}(D,{\bf{v}})$ be the complex vector space with basis the framed homogeneous tangles.  Notice
\begin{equation} \mathcal{S}(D,{\bf{v}})\leq \mathcal{F}(D,{\bf{v}}).\end{equation}

A {\bf framing} of a graph is an embedding of the graph in a compact oriented surface $\Sigma$ with boundary so that the graph is a deformation retract of $\Sigma$.  Framed homogeneous tangles correspond to isotopy classes of framed graphs relative to $\partial D^2 \times [0,1]$ of framed graphs embedded in $D^2\times [0,1]$ that always have the same side $\Sigma$ up where it touches $\partial D^2 \times [0,1]$. Computationally, we work on the level of diagrams, and use the fact that the value of the invariants should be unchanged by Reidemeister moves to derive equations that they must satisfy.

If $F:(D,{\bf{v}})\rightarrow (D',{\bf{w}})$ is an orientation preserving diffeomorphism so that $F({\bf{v}})={\bf{w}}$, and $F$ preserves the signs of the points, then $F$ gives rise to a linear map \begin{equation} F:\mathcal{S}(D,{\bf{v}})\rightarrow \mathcal{S}(D',{\bf{w}}).\end{equation}
Also,  $F$  extends to an orientation preserving diffeomorphism 
\begin{equation} \overline{F}:D\times [0,1]\rightarrow D'\times [0,1] \end{equation} by $\overline{F}(x,t)=(F(x),t)$. The map $\overline{F}$ induces a map 
\begin{equation} \overline{F}:\mathcal{F}(D,{\bf{v}})\rightarrow \mathcal{F}(D',{F(\bf{v}}))\end{equation}  You can understand the action of $\overline{F}$ on tangles from the action of $F$ on tangle diagrams.

\subsection{State Functions}

A choice $\{Z_D\}$  of $Z_D:\mathcal{F}(D,{\bf{v}})\rightarrow \mathcal{S}(D,{\bf{v}})$  for each marked oriented disk is {\bf topologically invariant} if for every orientation preserving diffeomorphism $F:D\rightarrow D'$, and  tangle diagram $T$, with boundary the marked set ${\bf{v}}$,
\begin{equation} F\circ Z_D(T)=Z_{D'}(\overline{F}(T)).\end{equation}

There is another type of symmetry  \cite{OY} that is assumed in the literature.  If $(D,{\bf{v}})$ is an oriented  disk along with a choice of signed points ${\bf{v}}$ in its boundary, let $(D,{\bf{v}}^*)$ denote the same disk and points, but with the signs of all the points changed.  If $ T\subset D$ is a tangle  diagram then let $T^*$ denote the tangle diagram having the same underlying set, but with the orientations of all the edges reversed.  There are linear maps,
\begin{equation} \mathcal{S}(D,{\bf{v}})\rightarrow  \mathcal{S}(D,{\bf{v}}^*)\ \mathrm{and} \  \mathcal{F}(D,{\bf{v}})\rightarrow  \mathcal{F}(D,{\bf{v}}^*)\end{equation} obtained by taking the linear extension of the map $T\rightarrow T^*$. We denote these by a super-scripted asterisk and call the operation the {\bf adjoint}. We say a choice of maps $\{Z_D\}$  {\bf preserves adjoints}  if for all  marked oriented disks $(D,{\bf{v}})$ and all framed homogeneous tangles  $T\in \mathcal{F}(D,{\bf{v}})$, 
\begin{equation} Z_D(T^*)=Z_D(T)^*.\end{equation}

Locality makes sense as before. Suppose that for every marked oriented disk $(D,{\bf{v}})$ we have chosen  a linear map \begin{equation} Z_D: \mathcal{F}(D,{\bf{v}})\rightarrow \mathcal{S}(D,{\bf{v}}).\end{equation}    If  $D'\subset D$ and $\partial D'$ is in general position  to the tangle diagram $\Gamma$, and  $\Gamma\cap \partial D'={\bf{w}}$ with decorations determined by the orientations of the edges of $F$ intersecting $\partial D'$.   Give $D'$ the orientation inherited from $D$. We can apply $Z_D'$ to $\Gamma\cap D'$ and then plug the result into $\Gamma \cap (D-D')$ to get an element of $\mathcal{F}(D,\bf{v})$ which is denoted  $Z_{D'}(\Gamma)$. The choice of functions $\{Z_D\}$  is {\bf local}  if for all $\Gamma$, $D'\subset D$ with $\Gamma $ in general position with  $\partial D'$,
\begin{equation} Z_D(\Gamma)=Z_D(Z_{D'}(\Gamma)).\end{equation}

\begin{definition} An {\bf $SU(3)$-state function} is a choice of linear maps $\{Z_D:\mathcal{F}(D,{\bf{v}})\rightarrow \mathcal{S}(D,{\bf{v}})\}$ for all marked oriented disks, that is local, topologically invariant, and preserves adjoints . \end{definition}

Finally, an $SU(3)$-state function  $\{Z_D\}$   is  {\bf separating} if all of the $Z_D$ are onto.
We do not know apriori that the inclusion map $\mathcal{S}(D,{\bf{v}})\rightarrow \mathcal{F}(D,{\bf{v}})$ is injective, so we have not declared $\mathcal{S}(D,{\bf{v}})\leq \mathcal{F}(D,{\bf{v}})$. This boils down to determining that if two webs differ by Reidemeister moves then they were isotopic to begin with.  If a state function exists, then the inclusion is injective, so eventually we can say that the state function is a projection.

\subsection{Twisting} The vector spaces $\mathcal{F}(D,{\bf{v}})$ and  $\mathcal{S}(D,{\bf{v}})$ have natural bases made up of equivalence classes of diagrams.  You can get different invariants by  changing what you mean by a diagram. That is you could interpret a diagram as meaning a scalar  multiple of itself. This will give rise to invariants that are different for individual diagrams, even though the state functions are similar as linear maps.  There are two kinds of twisting that make sense.  The first is twisting by vertices.  Since vertices ``represent'' invariant tensors, and a scalar multiple of an invariant tensor is invariant, you can multiply each diagram by a number based how many vertices it has, and still have a state function. The other is twisting by a function of crossings. In general twisting by a function of crossings does not make sense on the level of tangles.  This is because  the equations that $P$ must satisfy coming from Equation \ref{newr} are not homogeneous in the number of crossings. However, if you twist by a function of crossings, and restrict to diagrams that have no trivalent vertices, you get an invariant of tangles without trivalent vertices.

\begin{prop} Let $t\in \mathbb{C}-\{0\}$. For all $(D,{\bf{v}})$ let $P:\mathcal{F}(D,{\bf{v}})\rightarrow \mathcal{F}(D,{\bf{v}})$ be the linear map that sends any tangle diagram $T$ representing an element of $\mathcal{F}(D,{\bf{v}})$ to $t^rT$ where $r$ is the number of trivalent vertices of $T$. If $\{Z_D\}$ is a choice of maps that is a state function then so is  $\{P^{-1}\circ Z_D\circ P\}$. \end{prop}

\proof The number of trivalent vertices is preserved by diffeomorphisms, and Reidemeister moves. Since we do not distinguish between sources and sinks, $P$ commutes with taking the adjoint. \qed

We call $\{P^{-1}\circ Z_D \circ P\}$ the result of {\bf twisting $\{Z_D\}$ by $P$}.

\begin{prop} Let $c\in \mathbb{C}-\{0\}$. If $T$ is a tangle diagram let $p$ be the number of positive crossings of $T$ and $n$ be the number of negative crossings. Let $P(T)=c^{p-n}T$.  If the  map $P$ descends to a map of tangles, then $c^3=1$.  \end{prop}

\proof The only Reidemeister moves that are not homogeneous in the number of crossings come from Equation \ref{newr}.  We need to just look at one.
 \begin{equation}  P\left(\raisebox{-20pt}{\includegraphics{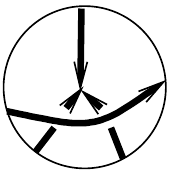}}\right)=P\left(\raisebox{-20pt}{\includegraphics{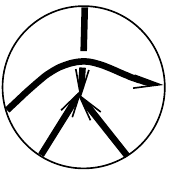}}\right)\end{equation}
Since there are two positive crossings on the left and one negative crossing on the right  for this equation to hold, $c^2=c^{-1}$. \qed

However, it is common to  twist tangle diagrams without vertices  according to the {\bf writhe}, to get an invariant of tangles  that is unchanged by the first Reidemeister move. 

\subsection{Uniqueness}

\begin{theorem} \label{bigenchilada} Every separating $SU(3)$-state function belongs to a specific  two parameter family of separating $SU(3)$-state functions parameterized by a pair of nonzero complex numbers $a$ and $y$.   \end{theorem}

\proof The proof of this theorem is a more elaborate version of the argument in the last section. We set up a system of equations and solve a system of equations derived from Reidemeister moves whose solutions correspond exactly to to separating $SU(3)$-state functions.  In the diagrams in this proof we assume that the disks lie in the plane and have inherited the standard orientation from $\mathbb{R}^2$.

In order to carry out the proof we need a complete list of non-elliptic webs with few boundary points. Specifically we need to understand $0$, $2$, $4$ and a special case of five boundary points. There are no webs with one boundary point, and the case of three boundary points isn't needed for the argument.

\begin{itemize}
\item $\mathcal{S}(D,\emptyset)$ is the vector space on the empty diagram.
  The fact that $Z_D$ preserves adjoints means that

  \begin{equation} Z_D\left(\raisebox{-20pt}{\includegraphics{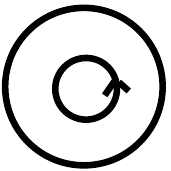}}\right)= Z\left(\raisebox{-20pt}{\includegraphics{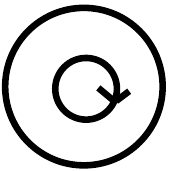}}\right)=\tau\emptyset\end{equation} for some $\tau \in \mathbb{C}$.

\item Suppose that ${\bf{v}}$ consists of two points with opposite signs. A basis for  $\mathcal{S}(D,{\bf{v}})$ is just,
  \begin{equation} \left\{\raisebox{-20pt}{\includegraphics{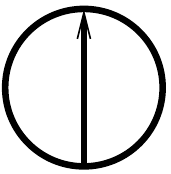}}\right\}. \end{equation}
  From this we know there is a complex number $\beta$ with 
 \begin{equation}\label{bubblez} Z_D\left(\raisebox{30pt}{\rotatebox{180}{\includegraphics{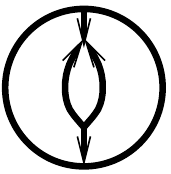}}}\right)=\beta\raisebox{-20pt}{\includegraphics{two.pdf}}.\end{equation}

\item Suppose that ${\bf{v}}$ consists of four points whose signs alternate as you traverse the circle then a basis for  $\mathcal{S}(D,{\bf{v}})$ is given by
  \begin{equation} \left\{\raisebox{-20pt}{\includegraphics{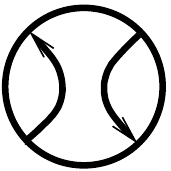}},\raisebox{-20pt}{\includegraphics{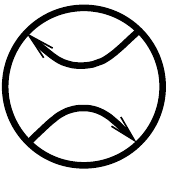}}\right\}. \end{equation}

  From this we conclude that there are complex numbers $a$ and $b$ so that
  \begin{equation}\label{squaresville} Z_D\left( \raisebox{-20pt}{\includegraphics{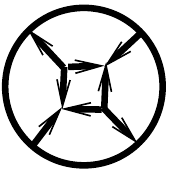}}\right)=a\raisebox{-20pt}{\includegraphics{alft41.pdf}}+b\raisebox{-20pt}{\includegraphics{alft42.pdf}}.\end{equation}

  Notice that if you rotate Equation \ref{squaresville} $\pi/2$ radians, the left hand side gets sent to its adjoint. The first diagram on the right gets sent to the adjoint of the second diagram on the right and vise verse. The fact that $Z_D$ preserves adjoints means that $a=b$ \cite{OY}.
  
\item Suppose that ${\bf{v}}$ consists of four points, two of which are decorated with $+$ and two of which are decorated with $-$ and the signs don't alternate, then a basis for  $\mathcal{S}(D,{\bf{v}})$ is given by
  \begin{equation} \left\{\raisebox{-20pt}{\includegraphics{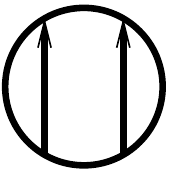}},\raisebox{-20pt}{\includegraphics{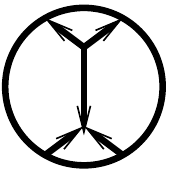}}\right\}. \end{equation}

  From this we conclude that there are $x,y,u,v\in \mathbb{C}$ so that

  \begin{equation}\label{posz}  Z_D\left(\raisebox{-20pt}{\includegraphics{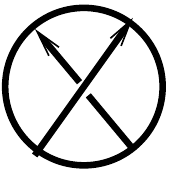}}\right)=x\raisebox{-20pt}{\includegraphics{fourn.pdf}}+y\raisebox{-20pt}{\includegraphics{fourv.pdf}},\end{equation}
 
  and

   \begin{equation}\label{negz}  Z_D\left(\raisebox{-20pt}{\includegraphics{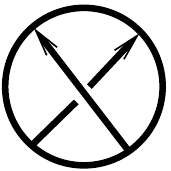}}\right)=u\raisebox{-20pt}{\includegraphics{fourn.pdf}}+v\raisebox{-20pt}{\includegraphics{fourv.pdf}},\end{equation}

  \item Suppose that ${\bf{v}}$ consists of five points, four of which carry the sign $+$ and one carries the sign $-$, then a basis is given by

 \begin{equation} \left\{\raisebox{-20pt}{\includegraphics{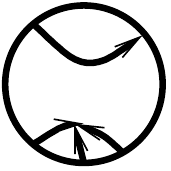}},\raisebox{-20pt}{\includegraphics{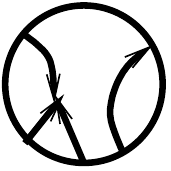}},\raisebox{-20pt}{\includegraphics{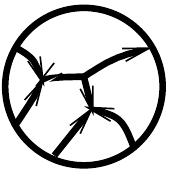}}\right\}. \end{equation}

  \end{itemize}

  Next we analyze the fact that the state function must be unchanged by three Reidemeister moves.  The first two are the second Reidemeister move with orientations added, and the third is a version of the new Reidemeister move from Equation \ref{newr}.

  \begin{itemize}
    \item From the second Reidemeister move we have
      \begin{equation}  Z_D\left(\raisebox{-20pt}{\includegraphics{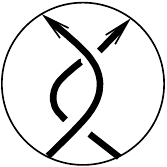}}\right)=Z_D\left(\raisebox{-20pt}{\includegraphics{fourn.pdf}}\right)=\raisebox{-20pt}{\includegraphics{fourn.pdf}}.\end{equation}
      Using Equations \ref{posz},\ref{negz} in small disks about the crossings and Equation \ref{bubblez} followed by  collecting terms we get,
      \begin{equation}\label{givesb} xu\raisebox{-20pt}{\includegraphics{fourn.pdf}}+(xv+uy+yv\beta)\raisebox{-20pt}{\includegraphics{fourv.pdf}}=\raisebox{-20pt}{\includegraphics{fourn.pdf}}\end{equation}
      From which we conclude $xu=1$ and $yv\beta=-xv-uy$.

      \item Now we resolve the second Reidemeister move where the two arcs are oriented oppositely to one another.

        \begin{equation}  Z_D\left(\raisebox{-20pt}{\includegraphics{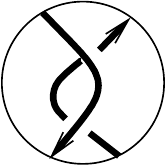}}\right)=Z_D\left(\raisebox{-20pt}{\includegraphics{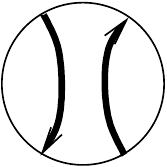}}\right)=\raisebox{-20pt}{\includegraphics{otherguy.pdf}}.\end{equation}

        Once again, applying Equations \ref{posz} and \ref{negz} to the left hand side in small disks about the crossings we get,

        \begin{equation}yv\raisebox{-20pt}{\includegraphics{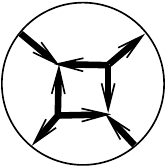}}+xv\raisebox{-20pt}{\includegraphics{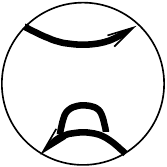}}+yu\raisebox{-20pt}{\includegraphics{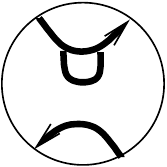}}+xu\raisebox{-20pt}{\includegraphics{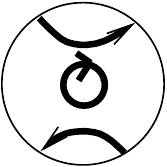}}  =\raisebox{-20pt}{\includegraphics{otherguy.pdf}}\end{equation}

        Using Equations \ref{bubblez} and \ref{squaresville}, remembering that $a=b$, we get

 \begin{equation}yva\raisebox{-20pt}{\includegraphics{otherguy.pdf}}+(xv\beta+uy\beta+xu\tau+yva)\raisebox{-20pt}{\includegraphics{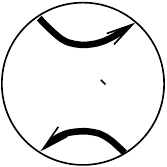}}  =\raisebox{-20pt}{\includegraphics{otherguy.pdf}}\end{equation}        

 From which we conclude $yva=1$ and $xu\tau+(xv+uy)\beta+yva=0$. This implies that $\tau=-(xv+uy)\beta-1$.

 \item Finally,

  \begin{equation}  Z_D\left(\raisebox{-20pt}{\includegraphics{newr2t.pdf}}\right)=Z_D\left(\raisebox{-20pt}{\includegraphics{newr2ta.pdf}}\right)\end{equation}
Resolving, using  our formulas, for crossings, the bubble and the square  and then collecting terms with respect to the basis for the disk with five points in the boundary four of which have the sign $+$ and one has the sign $-$ we get,

\begin{equation}(x^2+xy\beta+y^2a)\raisebox{-20pt}{\includegraphics{four2.pdf}}+(y^2a)\raisebox{-20pt}{\includegraphics{four1.pdf}}+xy\raisebox{-20pt}{\includegraphics{four3.pdf}}  =u\raisebox{-20pt}{\includegraphics{four1.pdf}}+v\raisebox{-20pt}{\includegraphics{four3.pdf}}.\end{equation}

From which we conclude,
\begin{equation} \label{redund}x^2+xy\beta+y^2a=0, \  y^2a=u,\ \mathrm{and} \ xy=v.\end{equation}

This means \begin{equation}\label{result}  u=y^2a, \ v=\frac{1}{ya}, \ \mathrm{and} \ x=\frac{1}{y^2a}. \end{equation}

\begin{itemize}
\item The first equation on the left of \ref{redund} is redundant, as it can be derived from the other equations and the formula for $\beta$ from Equation \ref{givesb}. 
\item  Notice that $x=v^2a$, $y=\frac{1}{va}$, and $u=\frac{1}{v^2a}$ so the system of Equations \ref{result}  is unchanged by simultaneously swapping $u$ and $x$ and $v$ and $y$. That means that all Reidemeister moves obtained  from the ones considered here by switching all the crossings yield no equations beyond the ones we have derived here.
\item Preservation of adjoints means if we reverse all the orientations on edges of the Reidemeister moves considered here we get the same set of equations. 
\item We  still need to consider Reidemeister moves from Equation \ref{newr} obtained by reversing the orientation on the edge that has no vertex. There is an additional symmetry  of the state function coming from rotating the diagrams in space along an axis lying in the plane which causes these equations to be  redundant. 

\item The third Reidemeister move is a consequence of the equations we have derived.
\end{itemize}
Since if you parameterize $x$,$u$ and $v$ by $a$ and $y$ as in Equation \ref{result}  all Reidemeister moves leave the value of $Z_D$ unchanged,  every choice of $a,y\in\mathbb{C}-\{0\}$ gives rise to a well defined choice of maps $Z_D:\mathcal{F}(D,{\bf{v}})\rightarrow \mathcal{S}(D,{\bf{v}})$, in the case of disks where the cardinality of ${\bf{v}}$ is $0,2$ or $4$.  These give rise to reduction rules like those in \cite{SW} (page 458). 

\begin{equation} S_1: \raisebox{-20pt}{\includegraphics{crossp.pdf}}\rightarrow \frac{1}{y^2a}\raisebox{-20pt}{\includegraphics{fourn.pdf}}+y\raisebox{-20pt}{\includegraphics{fourv.pdf}},\end{equation}
\begin{equation} S_2: \raisebox{-20pt}{\includegraphics{crossn.pdf}}\rightarrow y^2a\raisebox{-20pt}{\includegraphics{fourn.pdf}}+\frac{1}{ya}\raisebox{-20pt}{\includegraphics{fourv.pdf}},\end{equation}
\begin{equation}S_3 :\raisebox{-20pt}{\includegraphics{alft43.pdf}}\rightarrow a\raisebox{-20pt}{\includegraphics{alft41.pdf}}+a\raisebox{-20pt}{\includegraphics{alft42.pdf}},\end{equation}
\begin{equation} S_4:\raisebox{30pt}{\rotatebox{180}{\includegraphics{bubble.pdf}}}\rightarrow(-\frac{1}{y^3a}-y^3a^2)\raisebox{-20pt}{\includegraphics{two.pdf}}.\end{equation}
and
\begin{equation} S_5: \raisebox{-20pt}{\includegraphics{clockloop.pdf}}\rightarrow 
(a^3y^6+1+a^{-3}y^{-6}) \emptyset .\end{equation} 
These reduction rules are terminal and confluent, by an argument similar to the one in \cite{SW} for the $A_2$-spider.  

To define $Z_D$ for an arbitrary disk, start by applying $S_1$ and $S_2$ in small disks around the crossings. Follow this by iteratively surrounding any bubbles, squares or trivial components by small disks and applying $S_3,S_4$ and $S_5$ until the final result is a linear combination of non-elliptic webs.
\end{itemize}

\qed

\begin{rem} 
Let $U:D^2\rightarrow D^2$ be an orientation reversing diffeomorphism. We can extend $U$ to an orientation preserving map \begin{equation} \overline{U}:D^2\times [0,1]\rightarrow D^2\times [0,1] \end{equation} by $\overline{U}(x,t)=(U(x),1-t)$.  On the level of tangle diagrams this corresponds to applying $U$ and then switching all the crossings. Even though we did not require it, for any state function and any tangle \begin{equation} Z_D(\overline{U}(T))=U(Z_D(T)) .\end{equation} This is why half of the Reidemeister moves coming from Equation \ref{newr} are redundant.  On the categorical level this phenomenon was observed for all categories of diagrams with a single trivalent vertex in \cite{MPS}. \end{rem}

\begin{theorem} Any separating state function is twist equivalent to a state function with $a=1$. \end{theorem}

\proof Suppose that we have a state function with $a\neq 1$. Choose a fourth root $\zeta$  of $a$.  If $\Gamma$ is a non-elliptic web with $n$ trivalent vertices let $P(\Gamma)=(\frac{1}{\zeta})^n\Gamma$.  It is easy to check that $P^{-1}\circ Z_D\circ P$, satisfies the relations coming from $a'=1$, $y'=\sqrt{a}y$.  \qed

Since $P$ does not change basis vectors without a vertex, the link invariants coincide.

\begin{theorem} For each marked oriented disk $(D,{\bf{v}})$, the  map induced by inclusion,
\begin{equation} i:\mathcal{S}(D,{\bf{v}})\rightarrow \mathcal{F}(D,{\bf{v}})\end{equation}
is injective. After identifying $\mathcal{S}(D,{\bf{v}})$ with its image under inclusion, any state function is a projection. \end{theorem}

\proof Let $\{Z_D\}$ be a separating state function, which thanks to Theorem \ref{bigenchilada} we know exists.  Since $Z_D:\mathcal{F}(D,{\bf{v}})\rightarrow \mathcal{S}(D,{\bf{v}})$ is onto there is a finite collection of tangle diagrams $T_1,\ldots,T_n$ so that $\{Z_D(T_i)\}$ spans $\mathcal{S}(D,{\bf{v}})$.  For each tangle diagram $T_i$, place a small disk $D_c$ about each crossing, and apply $Z_{D_c}$, so that we have a linear combination of diagrams without crossings, whose image under $Z_D$ spans $\mathcal{S}(D,{\bf{v}})$. Next iterative place disks about trivial components, bigons, and squares, and apply the state function there. After finitely many steps, there is a collection of linear combinations of non-elliptic webs, whose image under $Z_D$ spans $\mathcal{S}(D,{\bf{v}})$.  This means that the restriction of $Z_D$ to $i\mathcal{S}(D,{\bf{v}}))$ maps onto $\mathcal{S}(D,{\bf{v}})$. Since $\mathcal{S}(D,{\bf{v}})$ is a finite dimensional vector space, this implies that the restriction of $Z_D\circ i$ is injective. This in turn implies that $i: \mathcal{S}(D,{\bf{v}})\rightarrow \mathcal{F}(D,{\bf{v}})$ is injective.

The proof  that $Z_D:\mathcal{F}(D,{\bf{v}})\rightarrow \mathcal{S}(D,{\bf{v}})$  is a projection is just a repeat of the proof of Lemma  \ref{project}.

\qed

Every $SU(3)$-state function is twist equivalent to the state function satisfying  the following skein relations for some $y\in \mathbb{C}-\{0\}$.

  \begin{equation}  Z_D\left(\raisebox{-20pt}{\includegraphics{crossp.pdf}}\right)=y^{-2}\raisebox{-20pt}{\includegraphics{fourn.pdf}}+y\raisebox{-20pt}{\includegraphics{fourv.pdf}},\end{equation}
 \begin{equation}  Z_D\left(\raisebox{-20pt}{\includegraphics{crossn.pdf}}\right)=y^2\raisebox{-20pt}{\includegraphics{fourn.pdf}}+y^{-1}\raisebox{-20pt}{\includegraphics{fourv.pdf}},\end{equation}
   \begin{equation} Z_D\left( \raisebox{-20pt}{\includegraphics{alft43.pdf}}\right)=\raisebox{-20pt}{\includegraphics{alft41.pdf}}+\raisebox{-20pt}{\includegraphics{alft42.pdf}}.\end{equation}
  \begin{equation} Z_D\left(\raisebox{30pt}{\rotatebox{180}{\includegraphics{bubble.pdf}}}\right)=(-y^3-y^{-3})\raisebox{-20pt}{\includegraphics{two.pdf}},\end{equation}
  and
  \begin{equation}Z_D\left(\raisebox{-20pt}{\includegraphics{clockloop.pdf}}\right)= Z_D\left(\raisebox{-20pt}{\includegraphics{countloop.pdf}}\right)=(y^6+1+y^{-6})\emptyset.\end{equation}
  
It satisfies the HOMFLY skein relation
\begin{equation}y^{-1} Z_D\left(\raisebox{-20pt}{\includegraphics{crossp.pdf}}\right)-yZ_D\left(\raisebox{-20pt}{\includegraphics{crossn.pdf}}\right)=(y^{-3}-y^3)
  \raisebox{-20pt}{\includegraphics{fourn.pdf}}\end{equation}

We call $Z_D$ the {\bf standard $SU(3)$-state function}. 
In the case of the $SU(3)$-link polynomial, the only other non-uniqueness is a choice of value for the empty link.

\subsection{Comparison with the Literature}

\begin{itemize}

\item {\em Kuperberg \cite{Ku} page 129.} Denote the state function from Kuperberg's invariant by $\{Z_D^K\}$.    For marked disk $(D,{\bf{v}})$,

\begin{equation} Z_D^K=Z_D|_{q^{1/6}}.\end{equation}
That is evaluate the standard state function at $y=q^{1/6}$.

Kuperberg's state function satisfies
  \begin{equation}q^{-1/6} Z^K_D\left(\raisebox{-20pt}{\includegraphics{crossp.pdf}}\right)-q^{1/6}Z^K_D\left(\raisebox{-20pt}{\includegraphics{crossn.pdf}}\right)=(y^{-1/2}-q^{1/2})
\raisebox{-20pt}{\includegraphics{fourn.pdf}}\end{equation}

\item {\em Ohtsuki and Yamada \cite{OY}, page 375.}  Denote the state function from Ohtsuki and Yamada's work by $\{Z_D^{OY}\}$.   For any  marked disk $(D,{\bf{v}})$,

\begin{equation} Z_D^{OY}=Z_D|_{-q^{-1}}.\end{equation}

. 
  It satisfies the  skein relation 
\begin{equation}q Z^{OY}_D\left(\raisebox{-20pt}{\includegraphics{crossp.pdf}}\right)-q^{-1}Z^{OY}_D\left(\raisebox{-20pt}{\includegraphics{crossn.pdf}}\right)=(q^{3}-q^{-3})
\raisebox{-20pt}{\includegraphics{fourn.pdf}}\end{equation}

\item {\em Sikora \cite{Si} page 869}  Sikora relates his invariants to Kuperberg's, which in turn relates it to the normalization here. Twist the standard state function by   $P(\Gamma)=y^{9v(\Gamma)/2}\Gamma$ where $v(\Gamma)$ is the number of trivalent vertices of $\Gamma$ and then substitute $y=-q^{-1/3}$ to get Sikora's invariant.
\begin{equation} Z_D^S=P^{-1}\circ Z_D\circ P|_{-q^{-1/3}}.\end{equation}

It satisfies the  skein relation
\begin{equation}q^{1/3} Z^S_D\left(\raisebox{-20pt}{\includegraphics{crossp.pdf}}\right)-q^{-1/3}Z^S_D\left(\raisebox{-20pt}{\includegraphics{crossn.pdf}}\right)=(q-q^{-1})\raisebox{-20pt}{\includegraphics{fourn.pdf}}\end{equation}

   \item {\em Khovanov \cite{Kh} page 1053} Khovanov is defining a framing independent version of the invariant. His invariant does not fit directly into the framework of state functions. Denote his invariant of a tangle diagram $\Gamma$ without trivalent vertices by $<\Gamma>_{Kh}$. 
   
  If $\Gamma$ is a diagram for a tangle without trivalent vertices, let $s(\Gamma)$ be the signed sum of the crossings, then
  
  \begin{equation} <\Gamma>_{Kh}=y^{8s(\Gamma)}Z_D(\Gamma)|_{q^{1/3}}.\end{equation} 
  \end{itemize}

  \end{document}